\documentclass[reqno,11pt]{amsart}
\usepackage{amsmath,amssymb,amscd,amsxtra,dsfont}
\usepackage{graphicx,pst-grad,pst-plot,pst-coil}
\usepackage{pstricks}
\usepackage{eucal}
\usepackage[margin=1.4in]{geometry}
\newtheorem{theorem}{Theorem}
\newtheorem{proposition}[theorem]{Proposition}
\newtheorem{cor}[theorem]{Corollary}
\newtheorem{lemma}[theorem]{Lemma}

\theoremstyle{remark}

\usepackage{mathptmx}
\usepackage[all]{xy}
\def\bfig{\vcenter\bgroup}
\def\efig{\egroup}

\renewenvironment{proof}{
{\medskip\par\noindent\it Proof.\ }}{\vskip 2ex\par}

\def \defi{\begin{definition}}
\def \edefi{\end{definition}}

\def \prop{\begin{proposition}}
\def \eprop{\end{proposition}}

\def \thm{\begin{theorem}}
\def \ethm{\end{theorem}}

\def \lem{\begin{lemma}}
\def \elem{\end{lemma}}

\def \cor{\begin{corollary}}
\def \ecor{\end{corollary}}

\newcommand{\ex}{\begin{example}}
\newcommand{\eex}{\end{example}}

\newenvironment{exer*}
  {\small\begin{exercise}}
  {\end{exercise}}

\newcommand{\probref}[1]{\textbf{\ref{#1}} } %%%Tu
\def \ex*{\begin{example*}}
\def \eex*{\end{example*}}

\newenvironment{remark*}{
{\medskip\noindent\it Remark.}}{\vskip 2ex\par}
\def \rem*{\begin{remark*}}
\def \erem*{\end{remark*}}

\newenvironment{claim*}{
{\medskip\noindent\it Claim.}}{\vskip 2ex\par}

\def \pf{\begin{proof}}
\def \epf{\end{proof}}

\def \enum{\begin{enumerate}}
\def \eenum{\end{enumerate}}

\numberwithin{equation}{section}
\numberwithin{figure}{section}

\newcommand{\C}{\mathbb{C}}

\newcommand{\GL}{\operatorname{GL}}
\newcommand{\Z}{\mathbb{Z}}
\newcommand{\ga}{\alpha}
\newcommand{\gb}{\beta}
\newcommand{\gc}{\gamma}

\newcommand{\calj}{\mathcal{J}}

\newcommand{\R}{\mathbb{R}}

\newcommand{\comp}{\mathrel{\scriptstyle\circ}}

\newcommand{\cinf}{{C^{\infty}}}

\newcommand{\Eu}{e^T}

\newcommand{\Fl}{\operatorname{F\ell}}

\newcommand{\hu}{\underline{h}}                             
\newcommand{\id}{\operatorname{id}}

\newcommand{\inv}{^{-1}}

\newcommand{\isomto}{\ \smash{\xrightarrow{\raisebox{-2ex}{\smash{\ensuremath{_{\sim} \mspace{-1mu}}}}}}\ }

\newcommand{\length}{\operatorname{length}}

\newcommand{\oh}{\bar{h}}

\newcommand{\Ortho}{\text{O}}

\newcommand{\pt}{\operatorname{pt}}
\newcommand{\Q}{\mathbb{Q}}

\newcommand{\Sym}{\operatorname{Sym}}
\newcommand{\Symp}{\operatorname{Sp}}

\newcommand{\term}[1]{\textbf{\textit{#1}}}

\newcommand{\tri}{\bigtriangleup}

\newcommand{\ty}{\tilde{y}}
\newcommand{\uh}{\underline{h}}
\newcommand{\Unit}{\text{U}}

\begin{document}
\title
{Computing the Gysin Map Using Fixed Points
}
\author{Loring W. Tu}
\address{Department of Mathematics\\
         Tufts University\\
         Medford, MA 02155, USA\\
         Tel.: 617 627 3262\\
         Fax: 617 627 3966}
\email{loring.tu@tufts.edu}
\urladdr{http://www.tufts.edu/\~{}ltu/}
\date{July 1, 2015; version 20}
\begin{abstract}
The Gysin map of a map between compact oriented manifolds
is the map in cohomology induced by the push-forward map in
homology.  In enumerative algebraic geometry, formulas for the Gysin map
of a flag bundle play a vital role.
These formulas are usually proven by algebraic or combinatorial means.
This article shows how the localization 
formula in equivariant cohomology provides a systematic 
method for calculating the Gysin homomorphism in the ordinary 
cohomology of a fiber 
bundle.  
As examples, we recover classical
pushforward formulas for generalized 
flag bundles.
Our method extends the classical formulas to fiber bundles
with equivariantly formal fibers.
%We also discuss a generalization, to compact Lie groups,
%of the Lagrange--Sylvester 
%symmetrizer and the Jacobi symmetrizer in interpolation 
%theory.
%\bigskip
%\noindent
%{\sc R\'esum\'e}.
%Cet article montre comment la formule 
%de localisation d'Atiyah--Bott--Berline--Vergne 
%fournit une m\'ethode syst\'ematique pour 
%calculer le morphisme de Gysin d'un fibr\'e 
%en cohomologie ordinaire.  Par exemple, 
%nous retrouvons les formules de pushforward 
%pour les fibr\'es de vari\'et\'es de drapeaux.  
%Nous donnons aussi une g\'en\'eralisation 
%aux groupes de Lie compacts des op\'erateurs 
%de sym\'etrisation de Lagrange--Sylvester et 
%de Jacobi dans la th\'eorie d'interpolation.
\end{abstract}
\maketitle

In enumerative algebraic geometry, to count the number of objects satisfying a
set of conditions, one method is to represent the objects satisfying
each condition by cycles in a parameter space $M$ and then to compute
the intersection of these cycles in $M$.
When the parameter space $M$ is a compact oriented manifold,
by Poincar\'e duality, the intersection of cycles can be calculated as
a product of classes in the rational cohomology ring.
Sometimes, a cycle $B$ in $M$ is the image $f(A)$ of a cycle $A$
in another compact oriented manifold $E$ 
under a map $f\colon E \to M$.
In this case the homology class $[B]$ of $B$ is the image
$f_*[A]$ of the homology class of $A$ under the
induced map $f_*\colon H_*(E) \to H_*(M)$ in homology,
and the Poincar\'e dual $\eta_B$ of $B$ is the image of the Poincar\'e
dual $\eta_A$ of $A$ under the map $H^*(E) \to H^*(M)$ in cohomology
corresponding to the induced map $f_*$ in homology.
This map in cohomology, also denoted by $f_*$, is called the
\term{Gysin map}, the \term{Gysin homomorphism}, 
or the \term{pushforward map} in cohomology.
It is defined by the commutative diagram
\[
\xymatrix{
H^k(E) \ar[r]^-{f_*} \ar[d]_-{{\rm P.D.}}^-{\simeq} & 
H^{k-(e-m)}(M) \ar[d]^-{{\rm P.D.}}_-{\simeq} \\
H_{e-k}(E) \ar[r]_-{f_*} & H_{e-k}(M),
}
\]
where $e$ and $m$ are the dimensions of $E$ and $M$ respectively and
the vertical maps are the Poincar\'e duality isomorphisms.
The calculation of the Gysin map
for various flag bundles plays an important 
role in enumerative algebraic geometry, for example in 
determining the cohomology classes of degeneracy loci 
(\cite{porteous}, \cite{jozefiak--lascoux--pragacz}, \cite{harris--tu}, 
and \cite[Ch.~14]{fulton}).
Other applications of the Gysin map, for example, to the
computation of Thom polynomials associated to Thom--Boardman singularities and to the computation of the dual cohomology classes of bundles of Schubert varieties, may be found in \cite{damon74}.

The case of a projective bundle associated to a vector 
bundle is classical  
\cite[Eq.~4.3, p.~318]{arbarello--cornalba--griffiths--harris}. 
Pushforward formulas for 
a Grassmann bundle
and for a complete flag bundle
are described in 
Pragacz \cite[Lem.\ 2.5 and 2.6]{pragacz}  
and Fulton and Pragacz \cite[Section 4.1]{fulton--pragacz}.  
For a connected reductive group $G$ with a
Borel subgroup $B$ and a parabolic subgroup $P$ containing $B$, 
Akyildiz and 
Carrell \cite{akyildiz--carrell} found a pushforward 
formula for the map $G/B \rightarrow G/P$.
For a nonsingular $G$-variety $X$ such that $X \rightarrow X/G$ is a principal 
$G$-bundle,
% What are the hypotheses on $G$?  Connected reductive algebraic group
Brion \cite{brion96} proved using representation theory 
a pushforward formula for the 
flag bundle $X/B \rightarrow X/P$.

The pushforward map for a fiber bundle
 makes sense more generally even 
if $E$ and $M$ are not manifolds (\cite[\S 8]{borel--hirzebruch} 
or \cite{chern}); for example, $E$ and $M$ may
be only CW-complexes, so long as the fiber $F$ is a 
compact oriented manifold.  For $G$ a compact connected Lie 
group, $T$ a maximal torus, and $BG$, $BT$ their respective classifying 
spaces, Borel and Hirzebruch found in 
\cite[Th.~20.3, p.~316]{borel--hirzebruch}
a pushforward formula for the universal 
bundle $BT \rightarrow BG$ with fiber $G/T$.

Unless otherwise specified, by cohomology we will mean singular cohomology 
with rational coefficients.
A $G$-space $F$ is said to be \term{equivariantly formal} 
if the canonical restriction map $H_G^*(F) \rightarrow H^*(F)$
from its equivariant 
cohomology to its ordinary cohomology 
is surjective.
The main result of this paper, Theorem~\ref{t:pushforward}, 
shows that 
if the fiber $F$ of a fiber bundle $E \to M$ is an equivariantly 
formal manifold and has finite-dimensional cohomology, then the 
Gysin map of the fiber bundle can be 
computed from the equivariant localization 
formula of Atiyah--Bott--Berline--Vergne for a torus action 
(\cite{atiyah--bott}, \cite{berline--vergne}).  
This provides a systematic method for calculating the 
Gysin map.  In particular, we recover all the pushforward 
formulas mentioned above, but in the differentiable category
instead of the algebraic category.

Equivariant formality describes a large class of 
$G$-manifolds whose equivariant cohomology behaves nicely 
\cite[note 5, pp.\ 185--186]{guillemin--sternberg}.  These 
manifolds include all those whose cohomology vanishes in 
odd degrees.  In particular, a homogeneous space
 $G/H$, where $G$ is a compact Lie group and $H$ a 
closed subgroup of maximal rank, is equivariantly formal.

In fact, the technique of this article applies more 
generally to fiber bundles whose fibers are not 
equivariantly formal.
Let $F_G$ be the 
homotopy quotient of the space $F$ by the group $G$ and 
$\pi\colon  F_G \rightarrow BG$ the associated fiber bundle with fiber 
$F$.  
For any fiber bundle $f\colon  E \rightarrow M$ with fiber $F$ 
and structure group $G$, there is a bundle map $(h,\hu)$ from the 
bundle $E\rightarrow M$ to the bundle $F_G \rightarrow BG$.  We say that a 
class  in $H^*(E)$ is an \term{equivariant 
fiber class} if it is in the image of $h^*\colon  H_G^*(F) \rightarrow 
H^*(E)$.  
In Theorem \ref{t:fiberclass} we compute the
pushforward of an equivariant 
fiber class of any fiber bundle $f\colon  E\rightarrow M$
such that the pullback $f^*\colon  H^*(M) \rightarrow H^*(E)$ is 
injective.

Using the residue symbol, Damon in \cite{damon73} computed the Gysin map for classical flag bundles, fiber bundles whose fibers are flag manifolds of the classical compact groups $\Ortho(n)$, $\Unit(n)$, and $\Symp (n)$.
Since these flag manifolds are equivariantly formal, our 
Theorem~\ref{t:formal} includes these cases, although in a different form.
In Section \ref{s:other}, we work out the case of $\Unit(n)$ as an example.

The pushforward formula in Theorem~\ref{t:pushforward} 
suggests a geometric interpretation and a generalization of
certain symmetrizing 
constructions in algebra.  To every compact connected Lie 
group $G$ of rank $n$ and closed subgroup $H$ of maximal 
rank, we associate a 
symmetrizing operator on the polynomial ring in $n$ 
variables.
When $G$ is the unitary group $\Unit(n)$ and $H$ is the 
parabolic subgroup $\Unit(k) \times \Unit(n-k)$ or the maximal torus
$\Unit(1) \times \dots \times \Unit(1)$ ($n$ times), this 
construction specializes to
the 
Lagrange--Sylvester symmetrizer and the Jacobi symmetrizer of 
interpolation theory respectively.

% Although our discussion is in terms of real 
% cohomology and takes place in the $\cinf$ category, using 
% the results of Edidin, Graham, and Totaro the same formalism 
% works for Chow groups in the algebraic category.  In this 
% way the pushforward formula (Theorem \ref{t:pushforward})
%  is also valid for the Chows
% groups of a fiber bundle of varieties.
% 
% The link between the pushforward in ordinary cohomology and 
% in equivariant cohomology is a two-way street.
% Indeed, one can transport 
% pushforward formulas in ordinary cohomology 
% to formulas in equivariant cohomology.
% In this way, for example, the pushforward formula for a 
% projective bundle can be interpreted as a formula for the 
% equivariant Chern numbers of a projective space under the 
% standard torus action.
% We will discuss this in a separate article.

This article computes the Gysin map of a fiber bundle with equivariantly formal fibers.  In the companion article \cite{pedroza--tu},
we compute the Gysin map of a $G$-equivariant map for a compact connected Lie group $G$.

It is a pleasure to acknowledge 
the support of the
Tufts Faculty Research Award Committee in 2007--2008
and
the hospitality of the 
Institut Henri Poincar\'{e}, the 
Universit\'{e} de Lille, and the Institut de 
Math\'{e}matiques de Jussieu.
I thank Michel Brion for explaining his work and for 
generously sharing some key ideas with me, and
Jeffrey D.\ Carlson for careful proofreading and valuable feedback.

\section{Universal Fiber Bundles}

We work in the continuous category until Section~\ref{s:pushforward}, 
at which point we will switch to the smooth category.
In this section, $G$ is a topological group and $f\colon E\rightarrow
M$ is a continuous
fiber bundle with fiber $F$ and structure group $G$.  This 
means $G$ acts on $F$ on the left and there is a principal $G$-bundle 
$P\rightarrow M$ such that $E\rightarrow M$ is the associated fiber bundle
$P \times_G F \rightarrow M$.  Recall that the 
\term{mixing space} $P \times_G F$
is the quotient of $P\times F$ by the diagonal action of 
$G$:
\begin{equation} \label{e:diagonal}
g\cdot (p,x) = (pg\inv, gx) \quad\text{for } (p,x) \in 
P\times F \text{ and } g \in G.
\end{equation}
We denote the equivalence class of $(p,x)$ by $[p,x]$.

Let $EG \rightarrow BG$ be the universal principal $G$-bundle.  One 
can form the associated fiber bundle $\pi\colon EG \times_G F \rightarrow 
BG$.  The space $F_G := EG \times_G F$ is called the 
\term{homotopy quotient} of $F$ by $G$, and 
its cohomology $H^*(F_G)$ is by definition the 
\term{equivariant cohomology} $H_G^*(F)$ of the $G$-space 
$F$.
The following lemma shows that the bundle $\pi\colon F_G
\to BG$ can serve as a universal fiber bundle with fiber $F$ and
structure group $G$.

\lem  \label{l:universal}
For any fiber bundle $f\colon E \rightarrow M$ with fiber $F$ and 
structure group $G$, there is a bundle map $(h, \hu)$ from 
$f\colon E \rightarrow M$ to $\pi\colon  F_G \rightarrow BG$ such that 
the bundle $E$ is isomorphic to the pullback bundle
$\hu^*(F_G)$.
\elem 

\pf 
The classifying map $\underline{h}$ of the principal bundle
$P\rightarrow M$ in the diagram
\[
\xymatrix{
P \ar[d] \ar[r] & EG \ar[d] \\
M \ar[r]_-{\underline{h}} & BG
}
\]
induces a map of fiber bundles

\begin{equation} \label{e:universal}
\begin{xy}
(2,10.5)*{E =};
(6,5)*{\xymatrix{
P\times_G F \ar[d]_-f \ar[r]^-h & 
EG \times_G F \ar[d]^-{\pi} \\
M \ar[r]_-{\underline{h}} & BG.}
};
(47.5,10)*{=F_G};
\end{xy}
\end{equation}
\vspace{.3in}

Recall that the \term{fiber product} over a space $B$ of two maps $\ga\colon M \to B$
and $\beta\colon N \to B$ is
\[
M \times_B N := \{ (x,y) \in M\times N \mid \ga(x) = \gb(y) \}
\]
and that the total space of the pullback to $M$ 
of a bundle $\beta\colon N \to B$
via a map $\ga\colon M \to B$ is $\ga^*N := M \times_B N$.
If $N$ is a right $G$-space for some topological group $G$,
then so is the fiber product $M \times_B N$, with
$(m,n)g = (m, ng)$ for $(m,n)\in M\times_B N$ and $g \in G$.

Since the principal bundle $P$ is isomorphic to the pullback $\hu^*(EG)$
of $EG$ by $\hu$, it is easily 
verified that $E$ is isomorphic to the pullback $\hu^*(F_G)$ of $F_G$ by $\hu$:
\begin{align*}
E &= P\times_G F \simeq \hu^*(EG) \times_G F = (M\times_{BG} 
EG) \times_G F \\
&\simeq M \times_{BG} (EG \times_G F) = M\times_{BG} F_G = 
\hu^*(F_G). 
\end{align*}
(In the computation above, the notation $\times_{BG}$ denotes the fiber product
and the notation $\times_G$ denotes the mixing construction,
and the isomorphism 
\[
(M\times_{BG} EG) \times_G F \simeq M \times_{BG} (EG \times_G F)
\] 
is given by
\[
\big[(m,e), x\big] \ \longleftrightarrow\  \big(m, [e,x]\big)
\]
for $m \in M$, $e \in EG$, and $x \in F$.)
\qed\epf

\section{Equivariant Formality}

Let $G$ be a topological group acting on a topological space $X$,
and $X_G$ the homotopy quotient of $X$ by $G$.
Since $X_G$ fibers over the classifying space $BG$ with fiber $X$,
there is an inclusion map $X \hookrightarrow X_G$ and
correspondingly a restriction homomorphism $H_G^*(X) \to H^*(X)$
in cohomology.
As stated in the introduction, the $G$-space $X$ is defined to be
\term{equivariantly formal} if this homomorphism $H_G^*(X)
\to H^*(X)$ is surjective; in this case, we also say that
every cohomology class on $X$ has an \term{equivariant
extension}.

The following proposition gives a large class of equivariantly formal spaces.

\prop \label{p:formal}
Let $G$ be a connected Lie group.
A $G$-space $X$ whose cohomology vanishes in odd degrees is equivariantly formal.
\eprop

\pf
By the homotopy exact sequence of the fiber bundle $EG \to BG$
with fiber $G$, the connectedness of $G$ implies that $BG$ is
simply connected.
Since $X_G \to BG$ is a fiber bundle with fiber $X$ over a simply
connected base space, the $E_2$-term of the spectral sequence of the fiber bundle $X_G \to BG$ is the tensor product
\[
E_2^{p,q} = H^p(BG) \otimes_{\Q} H^q(X)
\]
(see \cite[Th.\ 15.11]{bott--tu}).
Recall that the cohomology ring $H^*(BG)$ is a subring of a polynomial ring with 
even-degree generators \cite[\S 4]{tu}.
Thus, $H^p(BG) =0$ for all odd $p$.
Together with the hypothesis that
$H^q(X) =0$ for all odd $q$, this means that
the odd columns and odd rows of the $E_r$-terms will be 
zero for all $r$.
For $r$ even,
$d_r\colon E_r^{p,q} \to E_{r+1}^{p+r,q-r+1}$ 
changes the row parity (moves from an odd
row to an even row and vice versa);
for $r$ odd, $d_r$ changes the column parity.
Thus, all the differentials $d_r$ for $r \ge 2$ vanish,
so the spectral sequence degenerates at the $E_2$-term and
additively
\[
H_G^*(X) = H^*(X_G)  = E_{\infty} = E_2 = H^*(BG) \otimes_{\Q} H^*(X).
\]
This shows that $H_G^*(X) \to H^*(X)$ is surjective,
so $X$ is equivariantly formal;
in fact, for any $\ga \in H^*(X)$, the element $1 \otimes \ga \in H^*(BG) \otimes_{\Q} H^*(X)
= H_G^*(X)$ maps to $\ga$.
\qed\epf

\section{Fiber Bundles with Equivariantly Formal Fibers}
\label{s:formal}

In this section we compute the cohomology ring with rational coefficients of the total space of a fiber bundle with equivariantly
formal fibers. 
We assume tacitly that all spaces have a basepoint and that
all maps are basepoint-preserving.
By the fiber of a fiber bundle over a space $M$, 
we mean the fiber above the basepoint of $M$.
 
For any continuous fiber bundle $f\colon E \rightarrow M$ with fiber $F$ and group 
$G$, the diagram \eqref{e:universal} induces a commutative 
diagram of ring homomorphisms
\[
\xymatrix{
H^*(E) & H_G^*(F) \ar[l]_-{h^*} \\
H^*(M) \ar[u]^-{f^*} & H^*(BG). \ar[l]^-{\underline{h}^*} 
\ar[u]_-{\pi^*}
}
\]
Thus, both cohomology rings $H^*(M)$ and $H_G^*(F)$ are 
$H^*(BG)$-algebras, and we can form their tensor product 
over $H^*(BG)$.

\thm \label{t:formal}
Let $f\colon E \rightarrow M$ be a continuous fiber bundle with fiber $F$ and 
structure group $G$.  Suppose $F$ is equivariantly formal 
and its cohomology ring $H^*(F)$ is finite-dimensional.  Then 
\enum
\item[(i)] there is a ring isomorphism
\begin{align} \label{e:cohomology}
\varphi\colon  H^*(M) \otimes_{H^*(BG)} 
H_G^*(F) &\rightarrow H^*(E),\\
a \otimes b &\mapsto (f^*a)h^*b; \notag
\end{align}
\item[(ii)] the pullback map $f^*\colon H^*(M) \to H^*(E)$ is injective.
\eenum
\ethm 

\pf 
(i)  Because $E$ is isomorphic to the pullback $\underline{h}^*(F_G)$
of $F_G$, the map $h\colon  E \rightarrow F_G$ maps the fiber $F$ of $E$ isomorphically 
to the fiber of $F_G$.  Hence, the inclusion map of the 
fiber, $F \rightarrow F_G$, factors as
\[
F\rightarrow E \overset{h}{\to} F_G.
\]
This means that in cohomology, the restriction map $H^*(F_G) \rightarrow H^*(F)$ 
factors through $h^*$:
\[
H_G^*(F) \overset{h^*}{\to} H^*(E) \rightarrow H^*(F).
\]

Since the restriction $H_G^*(F) \rightarrow H^*(F)$ is surjective by 
the hypothesis of equivariant formality, there are classes $b_1, \ldots, 
b_r$ in $H_G^*(F)$ that restrict to a basis 
for $H^*(F)$. 
Then $h^*b_1, \ldots$, $h^*b_r$ are classes 
in $H^*(E)$ that restrict to a basis for $H^*(F)$.  By the 
Leray--Hirsch theorem (\cite[Th.\ 5.11 and Exercise 15.12]{bott--tu}
or \cite[Th.~4D.1, p.~432]{hatcher}), the cohomology $H^*(E)$ is a free 
$H^*(M)$-module with basis $h^*b_1, \ldots , h^*b_r$.

Next consider the fiber bundle $F_G \rightarrow BG$.  By the 
Leray--Hirsch theorem again, $H^*(F_G)$ is a free 
$H^*(BG)$-module of rank $r$ with basis $b_1, \ldots, 
b_r$.
It then follows that $H^*(M) \otimes_{H^*(BG)} H_G^*(F)$ is 
a free $H^*(M)$-module of rank $r$ with
basis $1\otimes b_1, \ldots, 1\otimes b_r$.  The ring 
homomorphism $\varphi$ in \eqref{e:cohomology} is a 
homomorphism of free $H^*(M)$-modules of the same rank.  
Moreover, $\varphi$ sends the basis $1\otimes b_1, \ldots, 
1\otimes b_r$ to the basis $h^*b_1, \ldots, h^*b_r$,
so $\varphi$ is an isomorphism.

\medskip
\noindent
(ii)  If $\{ a_i \}$ is a basis for $H^*(M)$, then $\{a_i \otimes 1 \}$ is part
of a basis for $H^*(M) \otimes_{H^*(BG)} H_G^*(F) \simeq H^*(E)$.
Hence, $f^*\colon H^*(M) \to H^*(E)$ is injective.
\qed\epf

Let $h\colon  E \rightarrow F_G$ be a map that covers the classifying 
map $\hu\colon  M \rightarrow BG$ of the fiber bundle $f\colon E \rightarrow M$.
By Theorem~\ref{t:formal}, a cohomology class in 
$H^*(E)$ is a finite linear combination of elements of the 
form $(f^*a)h^*b$, with $a \in H^*(M)$ and $b \in H_G^*(F)$.
Under the hypothesis that the fiber is a compact oriented manifold,
by the projection formula \cite[Prop.~6.15]{bott--tu},
\[
f_*\big((f^*a)h^*b\big) = af_*h^*b.
\]
Hence, to describe the pushforward $f_*\colon  H^*(E) \rightarrow 
H^*(M)$, it suffices to describe $f_*$ on the image of 
$h^*\colon  H_G^*(F) \rightarrow H^*(E)$.
Since $f^*\colon  H^*(M) \rightarrow H^*(E)$ is an injection,
it is in turn enough to give a formula for $f^*f_*h^*b$ for 
$b \in H_G^*(F)$.
This is what we will do in Section~\ref{s:pushforward}.

\section{The Relation between $G$-Equivariant Cohomology
and $T$-Equivariant Cohomology} \label{s:relation}

In the next two sections, let $G$ be a compact connected
Lie group acting on a manifold $F$
and $T$ a maximal torus in $G$.
Denote the normalizer of $T$ in $G$ by $N_G(T)$.
The \term{Weyl group} of $T$ in $G$ is the quotient group
$W :=N_G(T)/T$.
It is a finite reflection group.
The equivalence class in $W$ of an element $w \in N_G(T)$
should be denoted $[w]$, but in practice we use $w$ to denote both
an element of $N_G(T)$ and its class in $W$.
In a finite reflection group, every element $w$ is a product of reflections
and has a well-defined \term{length} $\length(w)$,
the minimal number of factors of $w$ when 
expressed as a product of reflections.
We define the \term{sign} of an element $w$ 
to be $(-1)^w := (-1)^{\length(w)}$.

The diagonal action of $G$ on $EG \times F$ in
\eqref{e:diagonal} may be written on the right as
\[
(e,x)g = (eg, g\inv x) \qquad\text{for } (e,x) \in EG \times F
\text{ and } g \in G.
\]
Since $EG=ET$, this action induces an action of the Weyl group
$W$ on the homotopy quotient
$F_T = ET \times_T F = (ET \times F)/T$:
\[
(e,x)T \cdot w = (e,x)w T \qquad\text{for } (e,x)T \in F_T \ \text{and}\ w \in W.
\]
(In general, if a Lie group $G$ containing a torus $T$ acts on the 
right on a space $Y$, then the Weyl group $W$ acts
on the right on the orbit space $Y/T$.)
It follows that there is an induced action of $W$ on 
$H_T^*(F)$.  
Again because $EG =ET$, there is a natural projection
$j\colon  F_T \rightarrow F_G$.  
As explained in \cite[Lemma 4]{tu}, since $j\colon F_T \to F_G$
is a fiber bundle with fiber $G/T$,
the induced map $j^*\colon  
H_G^*(F) \rightarrow H_T^*(F)$ identifies the $G$-equivariant cohomology 
$H_G^*(F)$ with the 
$W$-invariant elements of 
the $T$-equivariant cohomology $H_T^*(F)$.
In particular, $j^*$ is an injection.

For a torus $T$ of dimension $\ell$, the cohomology of its 
classifying space $BT$ is the polynomial ring 
\begin{equation} \label{e:cohombt}
H^*(BT) \simeq \Q [ u_1, \dots, u_{\ell}]
\end{equation}
(see \cite[\S 1]{tu}),
and $H^*(BG)$ is the subring of $W$-invariants:
\begin{equation} \label{e:cohombg}
H_G^*(\pt) = H^*(BG) \simeq \Q[u_1, \ldots, u_{\ell}]^W.
\end{equation}

\section{Pushforward Formula} \label{s:pushforward}

In this section, $G$ is a compact connected Lie group acting on a 
compact oriented manifold $F$, and $f\colon E\rightarrow M$ a $\cinf$ fiber bundle with 
fiber $F$ and structure group $G$.  
Let $T$ be a maximal 
torus in $G$.  The action of $G$ on the fiber $F$ restricts 
to an action of $T$ on $F$.  For simplicity we assume for 
now that 
the fixed point set $F^T$ of the $T$-action on $F$ consists of 
isolated fixed points.
(Note that $F^T$ is the fixed point set of $T$ on $F$,
while $F_T$ is the homotopy quotient of $F$ by $T$.)
For a fixed point $p \in F^T$, let $i_p\colon \{ p\} \rightarrow F$ be
the inclusion map and
\[
i_p^*\colon  H_T^*(F) \rightarrow H_T^*\big(\{ p \}\big) \simeq H^*(BT)
\]
the restriction map in equivariant cohomology. 
The normal bundle $\nu_p$ of $\{ p\}$ 
in $F$ is simply the tangent space $T_pF$ over the singleton space
 $\{ p\}$.
Since the torus $T$ acts on $T_p F$, the normal bundle
$\nu_p$ is a $T$-equivariant oriented vector bundle.
As such, it has an equivariant Euler class
$\Eu(\nu_p) \in H^*(BT)$,
which is simply the usual Euler class of the induced vector bundle
of homotopy quotients $(\nu_p)_T \to \{p\}_T = BT$.
At an isolated fixed point of a torus action, 
the equivariant Euler class $e^T(\nu_p)$ of the normal
bundle is nonzero and is therefore invertible in the fraction field of the
polynomial ring $H^*(BT)$ (see \cite[pp.~8--9]{atiyah--bott}).
For $b \in H_T^*(F)$, the fraction $(i_p^* b)/\Eu 
(\nu_p)$ is in the fraction field of $H^*(BT)$.

%Since $G$ acts on $ET \times F$ on the right,
%so does the subgroup $N_G(T)$.
%Thus, the Weyl group $W = N_G(T)/T$ acts on $(ET \times F)/T$
%on the right:  if $w \in N_G(T)$ and $xT \in F_T = (ET \times F)/T$, then
%$xTw = xwT$. 
%It follows that there is an induced action of $W$ on 
%$H_T^*(F)$.  Let $j\colon  F_T \rightarrow F_G$ be the natural projection 
%map.  By \cite[Prop.\ 1]{brion}, the induced map $j^*\colon  
%H_G^*(F) \rightarrow H_T^*(F)$ identifies the $G$-equivariant cohomology 
%$H_G^*(F)$ with the 
%$W$-invariant elements of 
%the $T$-equivariant cohomology $H_T^*(F)$.

  \lem  \label{l:injective}
Let $\pi\colon  F \rightarrow \pt$ be the constant map,
$\pi_G\colon F_G \to \pt_G=BG$ the induced map of homotopy quotients,
and $\pi^* = \pi_G^*\colon  
H_G^*(\pt) \rightarrow H_G^*(F)$ the induced map in $G$-equivariant 
cohomology.
If $F$ has a fixed point $p$, then $\pi^*$ is injective.
\elem 

\pf 
Let $i\colon  \pt \rightarrow F$ send the basepoint $\pt$ to the fixed 
point $p$.  Then $i$ is a $G$-equivariant map and $\pi 
\comp i = \id$.
It follows that $i^* \comp \pi^* = \id$ on $H_G^*(\pt)$.
Hence, $\pi^*$ is injective.
\qed\epf

Keeping the notations of Sections~\ref{s:formal} and \ref{s:relation},
we let $h\colon E \to F_G$ be a map that covers a classifying map
as in \eqref{e:universal} and $j\colon F_T \to F_G$ the natural
projection.

\thm \label{t:pushforward}
Let $f\colon  E \rightarrow M$ be a smooth fiber bundle with fiber $F$ and 
structure group $G$.  Let $T$ be a maximal torus in $G$.
%Suppose the pullback map $f^*\colon  H^*(M) \rightarrow 
%H^*(E)$ is injective.  
Suppose $F$ is a compact oriented equivariantly formal 
manifold and $T$ acts on $F$ with isolated fixed points.
Then for $b \in 
H_G^*(F)$, the rational expression
$\sum_{p \in F^T} (i_p^*j^*b)/\Eu(\nu_p)$ is in $H_G^*(\pt)$ 
and the pushforward map $f_*\colon H^*(E) \to H^*(M)$ is
completely specified by the formula
\begin{equation} \label{e:pushforward}
f^*f_*h^* b = h^* \pi^* \sum_{p \in F^T} \frac{i_p^* j^* 
b}{\Eu (\nu_p)},
\end{equation}
where the sum runs over all fixed points $p$ of the torus $T$ on $F$,
and $\pi^*:= \pi_G^*$ is the canonical map $H^*(BG) \to H_G^*(F)$.
(See diagram~\eqref{e:charmap} below for how the various maps
fit together.)
\ethm 

\rem* 
%For a torus $T$ of dimension $\ell$, the cohomology of its 
%classifying space $BT$ is a polynomial ring,
%\[
%H^*(BT) \simeq \Q [ u_1, \dots, u_{\ell}],
%\]
%and $H^*(BG)$ is the subring of $W$-invariants:
%\[
%H_G^*(\pt) = H^*(BG) \simeq \Q[u_1, \ldots, u_{\ell}]^W.
%\]
A priori, $(i_p^*j^* b)/\Eu(\nu_p)$ is a rational expression 
in $u_1, \dots, u_{\ell}$ (see \eqref{e:cohombt}).  However, it is part of the 
theorem that the sum $\sum_{p \in F^T} 
(i_p^*j^*b)/\Eu(\nu_p)$ is in fact a $W$-invariant 
polynomial in $u_1, \dots, u_{\ell}$, and hence is 
in $H_G^*(\pt)$.
\erem* 

\medskip
\noindent
\textit{Proof of Theorem \ref{t:pushforward}}.
For any $G$-space $X$, there is a natural projection $X_T 
\rightarrow X_G$ of homotopy quotients.  Hence, there is a 
commutative diagram
\begin{equation} \label{e:gtcomparison}
\bfig
\begin{xy}
(4.5,0)*{\xymatrix{
F_G \ar[d] & F_T \ar[l]_-j \ar[d]\\
{\ \ \pt_G} &**[r]{\pt_T} \ar[l]}
};
%(0,-13.5)*{BG =};
%(37, -13.5)*{=BT.}
(2,-13.5)*{BG =};
(35, -13.5)*{=BT.}
\end{xy}
\efig
\end{equation}
We append this commutative diagram to the commutative diagram 
arising from the classifying map of the fiber bundle $E \rightarrow 
M$:
\begin{equation} \label{e:classifying}
\bfig
\xymatrix{
E \ar[d]_-f \ar[r]^-h & F_G \ar[d]^-{\pi_{{}_G}} & F_T \ar[l]_-j \ar[d]^{\pi_{{}_T}} \\
M \ar[r]_-{\underline{h}} &BG & BT. \ar[l]
}
\efig
\end{equation}
By the push-pull formula 
(\cite[Prop.~8.3]{borel--hirzebruch} or \cite[Lem.~1.5]{chern}), this 
diagram induces a commutative diagram in cohomology
\begin{equation} \label{e:charmap}
\bfig
\xymatrix{
H^*(E) \ar[d]_-{f_*} & {H_G^*(F)\ } \ar[l]_-{h^*} \ar[d]^-{\pi_*} \ar@{>->}[r]^-{j^*} 
& H_T^*(F) \ar[d]^-{\pi_{{}_{T*}}}\\
H^*(M) & {H^*(BG)\ } \ar[l]^-{\underline{h}^*} \ar@{^{(}->}[r] & H^*(BT),
}
\efig
\end{equation}
where the two horizontal maps on the right are injections
by the discussion of Section~\ref{s:relation} and, to simplify the notation, we write $\pi_*$ for $\pi_{{}_{G*}}$.
Thus, for $b \in H_G^*(F)$,
\begin{equation} \label{e:pushpull}
f_* h^* b = \hu^* \pi_* b = \hu^*\pi_{{}_{T*}} j^*b.
\end{equation}

By the equivariant localization theorem for a torus action 
(\cite{atiyah--bott}, \cite{berline--vergne}),
\[
\pi_{{}_{T*}} j^* b = \sum_{p \in F^T} \frac{i_p^* j^* b}{\Eu 
(\nu_p)} \in H^*(BT).
\]
(The calculation is done in the fraction field of 
$H^*(BT)$, but the equivariant localization theorem guarantees
that the sum is in $H^*(BT)$.)
By the commutativity of the second square in \eqref{e:charmap},
$\pi_{{}_{T*}} j^* b \in H^*(BG)$.  
% Hence,
% \[
% \pi^* \pi_{{}_{T*}} j^* b = \pi^* \sum_{p \in F^T} \frac{i_p^* j^* b}
% {\Eu (\nu_p)} \in H_G^*(F).
% \]

Taking $f^*$ of both sides of \eqref{e:pushpull}, we 
obtain
\begin{align*}
f^*f_* h^* b &= f^* \hu^* \pi_{{}_{T*}} j^*b \\
&= h^* \pi^* \pi_{{}_{T*}} j^* b \\
&= h^* \pi^* \left( \sum_{p\in F^T} \frac{i_p^* j^* 
b}{\Eu(\nu_p)}\right). \tag*{\qed}
\end{align*}

\section{Generalizations of the Theorem}

On the total space $E$ of a fiber bundle $f\colon E \to M$ with
fiber $F$ and structure group $G$, there are two special types of 
cohomology classes:  (i) the pullback $f^*a$ of a class $a$
from the base, and (ii) the pullback $h^*b$ of a class $b$
from the universal bundle $F_G$ in the commutative
diagram~\eqref{e:universal}.
The first type is usually called a \term{basic class}.
For lack of a better term, we will call the second type an
\term{equivariant fiber class}.
According to Theorem~\ref{t:formal}, if the fiber $F$
of the fiber bundle is
equivariantly formal and has finite-dimensional cohomology,
then every cohomology class on $E$ is a finite linear combination
of products of basic classes with equivariant fiber classes.
Therefore, by the projection formula, to describe the pushforward map
$f_*\colon H_*(E) \to H_*(M)$, it suffices to describe
the pushforward $f_*(h^*b)$ of an equivariant fiber class $h^* b$.

While the hypothesis of equivariant formality is essential to 
describe completely the Gysin map in
Theorem~\ref{t:pushforward}, a closer examination reveals
that it is not needed for formula \eqref{e:pushforward} to hold.
In fact, formula \eqref{e:pushforward} holds for any smooth 
fiber bundle, with no hypotheses on the fiber.
In case $f^*\colon H^*(M) \to H^*(E)$ is injective, as in
Theorem~\ref{t:pushforward}, formula~\eqref{e:pushforward} 
determines $f_*(h^*b)$ and gives a pushforward formula
for the equivariant fiber class $h^*b$.
We state the conclusion of this discussion in the 
following theorem.

\thm \label{t:fiberclass}
Let $f\colon E \to M$ be a smooth fiber bundle with
fiber $F$ and structure group $G$.
Suppose a maximal torus $T$ in $G$ acts on $F$
with isolated fixed points.
Then for $b \in H_G^*(F)$,
\[
f^*f_*h^* b = h^* \pi^* \sum_{p \in F^T} \frac{i_p^* j^* 
b}{\Eu (\nu_p)}.
\]
In case the pullback $f^*\colon H^*(M) \to H^*(E)$ is
injective, this formula determines the pushforward $f_*(h^*b)$
of the equivariant fiber class $h^*b$ of $E$.
\ethm

If the fixed points of the $T$-action on the fiber $F$ are 
not isolated, Theorem \ref{t:pushforward} still holds 
provided one replaces the sum over the isolated fixed points 
with the sum of integrals over the components of the fixed 
point set,
\[
\sum_C \int_C \dfrac{i_C^* j^*b}{\Eu (\nu_C)},
\]
where $C$ runs over the 
components of $F^T$,
$i_C\colon C \to M$ is the inclusion map, and $\nu_C$ is the normal bundle to $C$ 
in $M$.
The Euler class $e^T(\nu_C)$ is nonzero \cite{atiyah--bott}, 
essentially because in the
normal direction $T$ has no fixed vectors, so that the
representation of $T$ on the normal space at any point has
no trivial summand.

Although the formula in Theorem~\ref{t:pushforward} looks
forbidding, it is actually quite computable.
In the rest of the paper, we will show how to derive various
pushforward formulas in the literature from Theorem~\ref{t:pushforward}.

\section{The Equivariant Cohomology of a Complete Flag Manifold} \label{s:complete}

In order to apply Theorem~\ref{t:pushforward} to a flag bundle,
we need to recall a few facts from \cite{tu} about the ordinary and equivariant cohomology of a complete flag manifold $G/T$, where
$G$ is a compact connected Lie group and $T$ a maximal torus in $G$.

A \term{character} of a torus $T$ is a multiplicative 
homomorphism $\gc\colon  T \rightarrow \C^{\times}$, 
where $\C^{\times}$ is the multiplicative
group of nonzero complex numbers.
If we identify $\C^{\times}$ with the general linear group $\GL(1,\C)$,
then a character is
a 1-dimensional complex representation 
of $T$.
Let $\hat{T}$ be the group of characters of $T$,
written additively:  if $\ga, \gb \in \hat{T}$ and $t \in 
T$, then we write
\[
t^{\ga} := \ga(t) \quad \text{and} \quad t^{\ga +\gb} := 
\ga(t)\gb(t).
\]

Suppose $X \rightarrow X/T$ is a principal $T$-bundle.  To each 
character $\gc$ of $T$, one associates a complex line bundle
$L(X/T, \gc)$ on $X/T$ by the mixing construction
\[
L(X/T, \gc) := X \times_{\gc} \C := (X \times \C)/T,
\]
where $T$ acts on $X \times \C$ by
\[
 (x, v) \cdot t = (x t, \gc(t\inv) v).
\]

Associated to a compact connected Lie group $G$
and a maximal torus $T$ in $G$ are two principal $T$-bundles:
the principal $T$-bundle $G \to G/T$ on $G/T$ 
and the universal $T$-bundle $ET \to BT$ on the classifying
space $BT$.
Thus, each character $\gamma\colon T \to \C^{\times}$ gives rise,
by the mixing construction, to a complex line bundle
\[
L_{\gamma} := L(G/T, \gamma) = G \times_{\gc} \C
\]
on $G/T$ and a complex line bundle 
\[
S_{\gamma} := L(BT, \gamma) = ET \times_{\gc} \C
\]
on $BT$.

The Weyl group $W$ of $T$ in $G$ acts on the character group $\hat{T}$ of 
$T$ by
\[
(w\cdot \gc)(t) = \gc(w\inv t w).
\]
If the Lie group $G$ acts on the right on a space $X$, then 
the Weyl group $W$ acts on the right on the orbit space $X/T$ by
\[
r_w (xT) = (xT) w = xwT.
\]
This action of $W$ on $X/T$ induces 
an action of $W$  
on the cohomology ring $H^*(X/T)$.
Moreover, for $w \in W$ and $\gamma \in \hat{T}$,
\[
r_w^* L_{\gamma} = L_{w\cdot \gamma}, \qquad  
r_w^* S_{\gamma} = S_{w\cdot \gamma}
\]
(see \cite[Prop.~1]{tu}).

Fix a basis $\chi_1, \dots, \chi_{\ell}$ for the 
character group $\hat{T}$, and let 
\[
y_i = c_1(L_{\chi_i}) \in H^2(G/T) \qquad \text{and}\qquad 
u_i = c_1(S_{\chi_i})\in H^2(BT)
\]
be the first Chern classes of the line bundles 
$L_{\chi_i}$ and $S_{\chi_i}$ on $G/T$ and on $BT$
respectively.
Then
\[
H^*(BT) = \Q [ u_1, \ldots, u_{\ell}].
\]
The Weyl group $W$ acts on the polynomial ring 
$\Q[ u_1, \ldots, u_{\ell}]$ by
\[
w\cdot u_i = w\cdot c_1(S_{\chi_i}) = c_1(S_{w\cdot\chi_i}).
\]
It acts on the polynomial ring $R:= \Q[ y_1, \ldots, y_{\ell}]$
in the same way.
The cohomology ring of $G/T$ is
\[
H^*(G/T) = \Q[ y_1, \ldots, y_{\ell}]/(R_+^W),
\]
where $(R_+^W)$ is the ideal generated by the homogeneous
$W$-invariant polynomials of positive degree in $R$
(see \cite[Th.~5]{tu}).
Since the cohomology of $G/T$ has only even-degree elements,
by Proposition~\ref{p:formal} the space $G/T$ is equivariantly formal
under the action of any connected Lie group.

% With $F=G/T$, the commutative diagram \eqref{e:gtcomparison} induces a
% commutative diagram
% \[
% \xymatrix{
% H_G^*(G/T) \ar[r]^-{j^*} & H_T^*(G/T) \\
% H^*(BG) \ar[u] \ar[r] & H^*(BT) \ar[u]_-{\pi_T^*}
% }
% \]
% in cohomology.
Consider the fiber bundle $(G/T)_T \to BT$ with fiber $G/T$.
Since $G/T$ is equivariantly formal and has finite-dimensional
cohomology, by Theorem~\ref{t:formal},
there is a ring isomorphism
\begin{align*}
\varphi\colon H^*(BT) \otimes_{H^*(BG)} H_G^*(G/T) &\isomto H_T^*(G/T),\\
a \otimes b &\mapsto (\pi_T^*a) j^*b.
\end{align*}
Now
\begin{equation} \label{e:gtg}
(G/T)_G = EG \times_G (G/T) \simeq (EG)/T = BT.
\end{equation}
Thus, 
\[
H_G^*(G/T) \simeq H^*(BT) = \Q[
u_1, \ldots, u_{\ell}].
\]
It is customary to denote $\varphi(u_i\otimes 1) = \pi_T^* (u_i) \in H_T^2(G/T)$
also by $u_i$, but we will write 
\[
\tilde{y}_i = \varphi(1 \otimes u_i) = j^* (u_i) \in H_T^*(G/T).
\]
Then the $T$-equivariant cohomology of $G/T$ may also be written
in the form
\[
H_T^*(G/T) \simeq \Q [ u_1, \ldots, u_{\ell}, \tilde{y}_1, \ldots, \tilde{y}_{\ell} ]/\calj,
\]
where $\calj$ is the ideal generated by $p(\tilde{y})- p(u)$ as $p$
runs over the invariant polynomials of positive degree in $\ell$ variables
\cite[Th.~11]{tu}.
Since $j^*\colon H_G^*(G/T) \to H_T^*(G/T)$ is a ring homomorphism,
we have
\begin{equation} \label{e:j*}
j^* b(u) = b(\tilde{y}_1, \ldots, \tilde{y}_{\ell}) =: b(\tilde{y}).
\end{equation}

%Recall that a $G$-equivariant vector bundle is a vector bundle
%with a $G$-action on both the total space and the base space
%such that the projection map is $G$-equivariant.
%For a character $\gamma \in \hat{T}$, the Lie group $G$
%acts on the line bundle $L_{\gamma} = G \times_{\gamma} \C$
%on the left by
%\[
%g \cdot [x,v] = [gx, v],
%\]
%making $L_{\gamma}$ into a $G$-equivariant line bundle over 
%$G/T$.
%Let $(L_\gamma)_T \to (G/T)_T$ be the induced complex
%line bundle on $(G/T)_T$ and let
%\[
%\tilde{y}_i = c_1\big( (L_{\chi_i})_T \big) \in H_T^2(G/T).
%\]
%Then the equivariant cohomology ring of $G/T$ under the action of $T$ on
%$G/T$ by left multiplication is
%\[
%H_T^*(G/T) \simeq \frac{\Q[u_1, \ldots, u_{\ell}, \tilde{y}_1, \ldots,
%\tilde{y}_{\ell}]}{\mathcal{J}},
%\]
%where $\mathcal{J}$ is the ideal in $\Q[u_1, \ldots, u_{\ell}, \tilde{y}_1,
%  \ldots,
%\tilde{y}_{\ell}]$ generated by $b(\tilde{y})-b(u)$ for all
%polynomials
%$b \in R_{+}^{W}$, the $W$-invariant polynomials of positive degree
%in $\ell$ variables.

The maximal torus $T$ acts on $G/T$ by left multiplication,
and the fixed point set is precisely the Weyl group $W = N_G(T)/T$.
At a fixed point $w\in W$, we have the following two formulas:
\begin{enumerate}
\item[(i)] (Restriction formula for $G/T$) \cite[Prop.\ 10]{tu} The restriction homomorphism
\[
i_w^*\colon H_T^*(G/T) \rightarrow H_T^*\big(\{w\}\big) \simeq H^*(BT)
\]
is given by
\[
i_w^* u_i = u_i, \qquad i_w^* \tilde{y}_i = w\cdot u_i.
\]
\item[(ii)] (Euler class formula) \cite[Prop.\ 13]{tu} The equivariant Euler class of the normal bundle $\nu_w$ at the fixed point $w \in W$ is
\[
e^T(\nu_w) = w \cdot \left( \prod_{\ga \in \tri^+} c_1(S_{\ga}) \right) 
=(-1)^w \prod_{\ga \in \tri^+} c_1(S_{\ga}) \in H^*(BT),
\]
where $\tri^+$ is a choice of positive roots of the adjoint
representation of $T$ on the complexified Lie algebra of $G$.
\end{enumerate}

\section{Complete Flag Bundles}

In this section $G$ is a compact connected Lie group
with maximal torus $T$, and $f\colon E \to M$ is a fiber
bundle with fiber $G/T$ and structure group $G$.
Let $X \to M$ be the associated principal $G$-bundle.
Then
\[
E = X \times_G (G/T) \simeq X/T
\]
and $ M \simeq X/G$, so the given bundle is isomorphic to 
$X/T \to X/G$.

With $F = G/T$ in the commutative diagram \eqref{e:classifying},
yielding
\begin{equation} \label{e:doublesquare}
\bfig
\xymatrix{
E \ar[d]_-f \ar[r]^-h & (G/T)_G \ar[d]^-{\pi_{{}_G}} & (G/T)_T \ar[l]_-j \ar[d]^{\pi_{{}_T}} \\
M \ar[r]_-{\underline{h}} &BG & BT, \ar[l]
}
\efig
\end{equation}
% we have
% \begin{equation} \label{e:gtg}
% (G/T)_G = EG \times_G (G/T) \simeq (EG)/T = (ET)/T = BT.
% \end{equation}
we see that the equivariant fiber classes on $E$ are of the form
$h^* b(u)$, where
\[
b(u) = b(u_1, \ldots, u_{\ell}) \in H_G^*(G/T) \simeq H^*(BT) =
\Q[u_1, \ldots, u_{\ell}].
\]

\thm  \label{t:gysingt}
For $b(u) \in H^*(BT) \simeq H_G^*(G/T)$, the 
pushforward of $h^*b(u)$ under $f$ is given by
\[
f^*f_*h^*b(u)
= h^* \sum_{w\in W} w \cdot \left(
\frac{b(u)}{\prod_{\ga\in \tri^+} c_1(S_{\ga})}\right)
=h^*\left(
\frac{\sum_{w\in W} (-1)^w w\cdot b(u)}
{\prod_{\ga \in \tri^+} c_1(S_{\ga})}  \right).
\]
\ethm 

\pf
%We first calculate the map $j^*\colon H_G^*(G/T) \to H_T^*(G/T)$.
%For any character $\gamma \in \hat{T}$, the homotopy quotients
%$(L_{\gamma})_G$ and $(L_{\gamma})_T$ fit into a
%commutative diagram
%\[
%\xymatrix{
%(L_{\gamma})_G \ar[d]  & (L_{\gamma})_T \ar[l] \ar[d] \\
%(G/T)_G  & (G/T)_T \ar[l].
%}
%\]
%It can be checked that $j^*(L_{\gamma})_G \simeq (L_{\gamma})_T$.
%Moreover,
%\[
%(L_{\gamma})_G = EG \times_G L_{\gamma} = EG \times_G (G \times_{\gamma} \C) \simeq EG \times_{\gamma} \C = S_{\gamma}.
%\]
%Hence,
%\begin{align*}
%j^* u_i &= j^* c_1(S_{\chi_i}) = j^* c_1 \big( (L_{\chi_i})_G\big)
%= c_1\big( j^*(L_{\chi_i})_G \big) \\
%&= c_1\big( (L_{\chi_i})_T \big) = \tilde{y}_i.
%\end{align*}
Because $G/T$ is equivariantly formal, Theorem~\ref{t:pushforward}
applies.
Under the action of $T$ on $G/T$ by left multiplication,
the fixed point set is the Weyl group $W = N_G(T)/T \subset G/T$,
so in Theorem \ref{t:pushforward}, $F^T = W$, a finite set, and
\begin{align*}
f^*f_*h^*b(u) &= h^*\pi^* \sum_{w\in W} 
\frac{i_w^* j^* b(u)}{e^T(\nu_w)} \\
&= h^*\pi^* \sum_{w\in W} 
\frac{i_w^* b(\tilde{y})}{e^T(\nu_w)}\qquad\qquad\qquad\qquad\qquad\qquad\text{\big(by \eqref{e:j*}\big)}\\
&= h^* \pi^* \sum_{w\in W} w \cdot \left(
\frac{b(u)}{\prod_{\ga\in \tri^+} c_1(S_{\ga})}\right)\\
&\qquad\qquad\qquad\qquad\text{(by the restriction and Euler class formulas)}\\
%&\qquad\big(\text{since } i_w^*j^*b(u) = i_w^* b(\ty) = w\cdot b(u)\big)\\
&=h^*\left(
\frac{\sum_{w\in W} (-1)^w w\cdot b(u)}
{\prod_{\ga \in \tri^+} c_1(S_{\ga})}  \right).
\end{align*}
In the last line we have omitted $\pi^*$ because it is injective (Lemma~\ref{l:injective}),
so that $H^*(BG)$ can be identified with a subring of $H_G^*(G/T)$.
\qed\epf

\section{The Characteristic Map}

If $T$ is a torus and $\hat{T}$ its character group, 
we let $\Sym(\hat{T})$ be the symmetric algebra of $\hat{T}$ over the field
$\Q$ of rational numbers; if $\chi_1, \ldots, \chi_{\ell}$ is a basis for
the character group $\hat{T}$, then
\[
\Sym(\hat{T}) = \Q [ \chi_1, \ldots, \chi_{\ell}].
\]
Associated to a principal $T$-bundle $X \to X/T$ is an algebra
homomorphism $c = c_{X/T}\colon \Sym(\hat{T})$ $\to H^*(X/T)$
called the \term{characteristic map} of $X/T$.

Each character $\gc \in \hat{T}$ gives rise to a complex line bundle
$L(X/T,\gamma) = X \times_{\gc} \C \to X/T$, as discussed 
in Section~\ref{s:complete}.
Define $c\colon \hat{T} \to H^2(X/T)$ by
\[
c(\gc) = \text{the first Chern class}\  c_1\big( L(X/T, \gc)\big) \in H^2(X/T).
\]
This map can be checked to be a homomorphism of abelian groups
\cite[Section 1]{tu}.
The extension of this group homomorphism to an algebra homomorphism
$c\colon \Sym(\hat{T}) \to H^*(X/T)$ is the characteristic map of $X/T$.

%In order to apply Theorem~\ref{t:pushforward} to a fiber bundle
%with a homogeneous space as fiber,
%we first recall a few facts about the 
%characteristic map and the actions of the Weyl group
%from \cite[\S 2]{tu}.
%
%
%Define $c\colon  \hat{T} \rightarrow H^2(X/T)$ by
%\[
%c(\gc) = c_1(L(X/T, \gc)) \in H^2 (X/T).
%\]
%Let $\Sym(\hat{T})$ be the symmetric algebra of $\hat{T}$ 
%over the field $\Q$ of rational numbers.  
%The group homomorphism $c\colon  \hat{T} \rightarrow H^2(X/T)$ extends to 
%an algebra homomorphism $c\colon  \Sym(\hat{T}) \rightarrow H^*(X/T)$, called the 
%\term{characteristic map} of the bundle $X \to X/T$, also written 
%$c_{X/T}$.  

The associated line bundles satisfy the following 
functorial property.  

\lem  \label{l:pullback}
Let $(\bar{h}, h)$ be a $T$-equivariant 
bundle map from $X \rightarrow X/T$ to 
$Y\rightarrow Y/T$.
For each character $\gc$ of $T$, the map $h$ pulls the 
bundle $L(Y/T, \gc)$ back to $L(X/T, \gc)$:
\[
h^*L(Y/T, \gc) \simeq L(X/T, \gc).
\]
\elem 

\pf 
An element of $h^*L(Y/T,\gc)$ is an ordered pair $\big(xT, [y,v']_T\big)$ such that
$yT = h(xT)= \bar{h}(x)T$.  Hence, $y=\bar{h}(x)t$ for some $t \in T$
and
\[
[y, v']_T = [\bar{h}(x) t, v']_T = [\bar{h}(x), tv'] = [\bar{h}(x), v]_T,
\]
where we set $v = tv' = \gamma(t) v'$.

The map $\phi\colon  h^*L(Y/T, \gc) \rightarrow L(X/T, \gc)$,
\[
\big(xT, [\bar{h}(x), v]\big) \mapsto [x,v],
\]
is a well-defined bundle map and has an obvious inverse.
\qed\epf

It follows from this lemma that the characteristic map also 
satisfies a functorial property.

\lem  \label{l:commute}
Under the hypotheses above, the diagram
\begin{equation} \label{e:diagram}
\bfig
\xymatrix{
& \Sym(\hat{T}) \ar[dl]_-{c_{X/T}} \ar[dr]^-{c_{Y/T}} & \\
H^*(X/T) && H^*(Y/T) \ar[ll]^-{h^*}
}
\efig
\end{equation}
is commutative.
\elem 

\pf
For $\gamma \in \hat{T}$,
\[
h^*\big(c_{Y/T}(\gamma)\big) = h^*\Big( c_1\big(L(Y/T,\gamma)\big)\Big) 
= c_1\Big(h^*\big(L(Y/T,\gamma)\big)\Big)
= c_1\big(L(X/T, \gamma)\big) = c_{X/T}(\gamma).
\]
Since $h^* \comp c_{Y/T}$ and $c_{X/T}$ are both algebra homomorphisms
and $\Sym(\hat{T})$ is generated by elements of $\hat{T}$, the lemma follows.
\qed\epf

Suppose a compact Lie group $G$ with maximal torus $T$ acts on the
right on two spaces $X$ and $Y$ in such a way that $X \to X/G$ and $Y \to Y/G$
are principal $G$-bundles.
Then $X \to X/T$ and $Y\to Y/T$ are principal $T$-bundles, and
the Weyl group $W = N_G(T)/T$ acts on $\hat{T}$, $X/T$, and $Y/T$,
thus inducing actions on $\Sym \hat{T}$, $H^*(X/T)$, and $H^*(Y/T)$.
By \cite[Cor.~2]{tu}, the characteristic maps $c_{X/T}$ and $c_{Y/T}$ are
$W$-homomorphisms.
If $\bar{h}\colon X \to Y$ is a $G$-equivariant map,
$h\colon X/T \to Y/T$ is the induced map,
and $r_w$ and $r'_w$ are right actions of $w\in W$ on $X/T$ and $Y/T$ respectively,
then $h\comp r_w = r'_w \comp h$, so
the induced map $h^*\colon H^*(Y/T) \to H^*(X/T)$
in cohomology is also a $W$-homomorphism.\label{p:whomomorphism}
Thus, all three maps in the commutative diagram~\eqref{e:diagram} are $W$-homomorphisms.

\lem \label{l:pullbackquot}
Suppose a group $G$ containing a subgroup $T$ acts on the right
on two spaces $X$ and $Y$, and $\bar{h}\colon X \to Y$ is a $G$-equivariant map.
If $\underline{h}\colon X/G \to Y/G$ is the induced map of quotients,
then the pullback by $\hu$ commutes with the quotient by $T$:
\[
(\hu^*Y)/T = \hu^*(Y/T).
\]
\elem

\pf
By inserting quotients by $T$ in the pullback diagram,
we have a commutative diagram
\[
\xymatrix{
\hu^* Y \ar[r] \ar[d] & Y \ar[d]\\
(\hu^*Y)/T \ar[r] \ar[d] & Y/T \ar[d]\\
X/G \ar[r]_-{\hu}  & Y/G.
}
\]
By the definition of pullback,
\[
\hu^* Y = \{ (xG, y) \in X/G \times Y \mid \bar{h}(x) G = yG \}.
\]
Hence,
\[
(\hu^*Y)/T = \{ (xG, y) T = (xG, yT) \in X/G \times Y/T \mid  \bar{h}(x) G = yG \}.
\]
On the other hand,
\[
\hu^*(Y/T) = \{ (xG, yT) \in X/G \times Y/T \mid  \bar{h}(x) G = yG \}.
\]
Thus,
\[
(\hu^*Y)/T = \hu^*(Y/T). \tag*{\qed}
\]
\epf

Now let $G$ be a compact Lie group with maximal torus $T$
and $X \to X/G$ a principal $G$-bundle.
Let $\hu\colon X/G \to BG$ be the classifying map of $X \to X/G$,
so that there is a commutative diagram
\[
\xymatrix{
X \ar[r]^-{\bar{h}} \ar[d] & EG \ar[d]\\
X/G \ar[r]_-{\hu} & BG
}
\]
with $X \simeq \hu^*(EG)$.
Let $h\colon X/T \to (EG)/T$ be the map of quotients induced
from $\bar{h}$.
By Lemma~\ref{l:pullbackquot} and \eqref{e:gtg},
\[
X/T \simeq (\hu^* EG)/T = \hu^*(EG/T) = \hu^*(BT) \simeq \hu^*\big( (G/T)_G\big).
\]
We therefore have the commutative diagram
\[
\begin{xy}
(6,5)*{\xymatrix{
X/T \ar[r]^-h \ar[d]_-f & (EG)/T  \ar[d]\\
X/G \ar[r]_-{\hu} & BG.}
};
(53,10)*{=BT \simeq (G/T)_G};
\end{xy}
\]
\vspace{.3in}

In Theorem~\ref{t:gysingt}, let $b(u)$ be the characteristic class
$c_{ET/T}(\gc) = c_1(S_{\gc}) \in H^*(BT)$ 
for some $\gc \in \Sym(\hat{T})$.  By Lemma \ref{l:commute},
in which we take $Y=EG=ET$,
\[
h^*c_{ET/T}(\gc) = c_{X/T}(\gc).
\]
Because $h^*$ and $c_{X/T}$ commute with the action of the 
Weyl group, Theorem~\ref{t:gysingt} 
becomes
\begin{align*}
f^*f_* c_{X/T} (\gc) &= 
\frac{\sum_{w\in W} (-1)^w w\cdot h^* c_{ET/T}(\gc)}
{\prod_{\ga \in \tri^+} h^* c_{ET/T}(\ga)}
= \frac{\sum_{w\in W} (-1)^w w\cdot  c_{X/T}(\gc)}
{\prod_{\ga \in \tri^+} c_{X/T}(\ga)}
 \\
&=\frac{c_{X/T} \big(\sum_{w\in W} (-1)^w w\cdot \gc \big)}
{c_{X/T} \big({\prod_{\ga \in \tri^+} \ga}\big)},
\end{align*}
which agrees with Brion's pushforward formula for a complete flag 
bundle \cite[Prop.\ 1.1]{brion96}, with the difference that our formula is in the differentiable category with $G$ a compact connected Lie group, while Brion's formula is in the algebraic category with $G$ a reductive connected algebraic group.

\section{Partial Flag Bundles}

Keeping the notations of the preceding two sections, let $H$ be a 
closed subgroup of the compact connected Lie group
$G$ containing the maximal torus $T$.
The map $f\colon X/H \rightarrow X/G$ is a fiber bundle with fiber $G/H$ 
and structure group $G$, with $G$ acting on $G/H$ by left 
multiplication.  
Since $G/H$ has cohomology only in even degrees \cite[Th.~6]{tu},
it is equivariantly formal,
so Theorem \ref{t:pushforward} suffices to describe the 
Gysin map of $f$.

Denote by $W_H$ and $W_G$ the Weyl groups of $T$ in $H$ and 
in $G$ respectively.
By Lemma~\ref{l:universal} there is a bundle map $(h,\hu)$ 
from the fiber bundle $f\colon  X/H \rightarrow X/G$ to the fiber 
bundle $\pi\colon  (G/H)_G \rightarrow BG$.  The cohomology of 
$(G/H)_G = (EG) \times_G G/H = EG/H =BH$ is
\[
H_G^*(G/H) = H^*(BH) = H^*(BT)^{W_H} = \Q [u_1, \dots, 
u_{\ell}]^{W_H},
\]
the ring of $W_H$-invariant real polynomials in $u_1, \ldots, u_{\ell}$
(see \eqref{e:cohombg}).
Choose a set $\tri^+(H)$ of positive roots of $H$ and a 
set $\tri^+$ of positive roots of $G$ containing 
$\tri^+(H)$.

\thm  \label{t:gysingh}
For $b(u) \in H^*(BH)$, the pushforward of the equivariant fiber class $h^*b(u) \in 
H^*(X/H)$ under $f$ is given by
\[
f^*f_* h^* b(u) = h^* \sum_{w \in W_G/W_H} w\cdot 
\left( \frac{b(u)}
{\prod_{\ga \in \tri^+ - \tri^+(H)} c_1(S_{\ga})} \right).
\]
\ethm 

\pf 
By \cite[Prop.\ 14, Th.\ 11(ii), Th.\ 19]{tu} we have the following 
facts concerning the equivariant cohomology of $G/H$:
\begin{enumerate}
\item[(i)] The fixed point set of the action of $T$ on $G/H$
by left multiplication is 
\[
W_G/W_H = N_G(T)/N_H(T) = N_G(T)/\big(N_G(T)\cap H\big)
\subset G/H.
\]
\item[(ii)] The $T$-equivariant cohomology of $G/H$ is
\[
H_T^*(G/H) = \big(\Q [u_1, \dots, u_{\ell}] \otimes \Q[\ty_1, 
\dots, \ty_{\ell}]^{W_H}\big)/\mathcal{J},
\]
where $\mathcal{J}$ is the ideal generated by $p(\ty)-p(u)$ as $p$ 
ranges over all $W_G$-homogeneous polynomials of positive 
degree in $\ell$ variables.
\item[(iii)] (Restriction formula for $G/H$) If $i_w\colon \{w\} \hookrightarrow G/H$
is the inclusion map of a fixed point $w\in W_G/W_H$, then the
restriction homomorphism
\[
i_w^*\colon H_T^*(G/H) \to H_T\big(\{w\}\big) \simeq H^*(BT)
\]
in equivariant cohomology is given by
\[
i_w^* u_i = u_i, \qquad i_w^* f(\ty) = w\cdot f(u)
\]
for any $W_H$-invariant polynomial $f(\ty) \in \Q[\ty_1, \ldots, \ty_{\ell}]^{W_H}$.
\item[(iv)] (Euler class formula) The equivariant Euler class of the normal 
bundle $\nu_w$ at a fixed point $w\in W_G/W_H$ is
\[
e^T(\nu_w) = w\cdot \left(
\prod_{\ga \in \tri^+ - \tri^+(H)} c_1(S_{\ga}) \right).
\]
\end{enumerate}

By plugging these facts into Theorem~\ref{t:pushforward}, the theorem follows as in the proof 
of Theorem~\ref{t:gysingt}.
\qed\epf

%\rem* 
%There are no characteristic maps $c_{BH}\colon  \Sym \hat{T} \rightarrow 
%H^*(BH)$ or $c_{X/H}\colon  \Sym \hat{T} \rightarrow H^*(X/H)$, because 
%$BH$ and $X/H$ are not of the form $Y/T$.
%However, we have the following proposition.
%\erem* 
%
%\prop 
%If $\gc \in (\Sym \hat{T})^{W_H}$ is a $W_H$-invariant
%polynomial in the characters of $T$, then $c_{X/T}(\gc) \in 
%H^*(X/T)$ actually lands inside $H^*(X/H)$ $= H^*(X/T)^{W_H} \subset H^*(X/T)$.
%\eprop 
%
%\pf 
%By the hypothesis that $\gamma \in \Sym(\hat{T})^{W_H}$,
%we have $w \cdot \gamma = \gamma$ for all $w \in W_H$.
%Because the action of the Weyl group $W_H$ commutes with 
%the characteristic map $c_{X/T}$, for all $w\in W_H$,
%\[
%w \cdot c_{X/T}(\gc) = c_{X/T}(w\cdot \gc) = c_{X/T}(\gc).
%\]
%Hence, $c_{X/T}(\gamma)$ is a $W_H$-invariant element of $H^*(X/T)$:
%\[
%c_{X/T}(\gc) \in H^*(X/T)^{W_H} = H^*(X/H). \tag*{\qed}
%\]
%\epf
%
%
%For $\gc \in (\Sym \hat{T})^{W_H}$ and $b(u)= c_{ET/T}(\gc)\in H^*(BT)^{W_H}$,
%the commutative diagram \eqref{e:diagram} with $Y=ET$ gives
%\[
%h^*b(u) = h^*c_{ET/T}(\gc) = h^*c_{ET/T}(\gamma) = c_{X/T}(\gc).
%\]
%Since $h^*$ and $c_{X/T}$ both commute with the action of the Weyl group,
%Theorem~\ref{t:gysingh} becomes
%\begin{align*}
%f^*f_* c_{X/T}(\gc) &= \sum_{w\in W_G/W_H} w \left(
%\frac{h^* c_{ET/T}(\gc)}
%{\prod_{\ga \in \tri^+ - \tri^+(H)} h^*c_{ET/T}(\ga)} \right)\\
%&= c_{X/T} \left(
%\sum_{w\in W_G/W_H} w \left(
%\frac{\gc}
%{\prod_{\ga \in \tri^+ - \tri^+(H)} \ga} \right)  \right),
%\end{align*}
%which agrees with Brion's pushforward formula for a partial 
%flag bundle \cite[Prop.~2.1]{brion96}.

\section{Other Pushforward Formulas} \label{s:other}

In this section we show that the Borel--Hirzebruch 
formula \cite{borel--hirzebruch} may be derived in the same
manner as Theorem~\ref{t:pushforward} and that the formulas of
Fulton--Pragacz \cite{fulton--pragacz} for a complete flag bundle
and Pragacz \cite{pragacz} for a Grassmann bundle
are consequences of Theorem~\ref{t:pushforward}.

\subsection{The Borel--Hirzebruch Formula}

As before, $G$ is a compact connected Lie group with maximal 
torus $T$.  Let $EG \rightarrow BG$ and $ET \rightarrow BT$ be the universal 
principal $G$-bundle and $T$-bundle respectively.
Since $BT=(ET)/T= (EG)/T$ and $BG = (EG)/G$, the natural 
projection $\pi\colon  BT \rightarrow BG$ is a fiber bundle with fiber 
$G/T$. From Theorem~\ref{t:pushforward} we will deduce a formula of Borel and 
Hirzebruch for the Gysin map of $BT \rightarrow BG$.
Although the Borel--Hirzebruch formula concerns a fiber bundle
with a homogeneous space $G/T$ as fiber, it is not a special case
of the formulas of Akyildiz--Carrell \cite{akyildiz--carrell} or Brion \cite{brion96},
because $BT$ and $BG$ are infinite-dimensional.
It is, however, amenable to our method, because $BT$ and $BG$ are homotopy quotients of finite-dimensional manifolds by the group $G$.

Let $W$ be 
the Weyl group of $T$ in $G$.
Let $\ga_1, \dots, \ga_m$ be a choice of positive roots for $T$ in $G$, 
and write $a_i=c_1(S_{\ga_i}) \in H^2(BT)$ 
for their images under the characteristic map. 
 
\thm[\cite{borel--hirzebruch}, Th.~20.3, p.~316]
For $x \in H^*(BT)$, the pushforward under $\pi_*$ is
\[
\pi_* x = \dfrac{\sum_{w\in W} (-1)^w w\cdot x}{a_1 \cdots 
a_m}.
\]
\ethm 

\pf 
If we represent $BT$ as the homotopy quotient $(G/T)_G$
and $BG$ as the homotopy quotient $(\pt)_G$, then there 
is a commutative diagram
%\[
%\xymatrix{
%**[l]BT = (G/T)_G \ar[d]_-{\pi} & (G/T)_T \ar[l]_-j \ar[d]^{\pi_{{}_T}}\\
%**[l]BG = (\pt)_G &  **[r](\pt)_T = BT.  \ar[l]
%}
%\]

\bigskip

\hspace{1.5in}\begin{xy}
(0,10.5)*{BT =};
(6.5,5)*{\xymatrix{
(G/T)_G \ar[d]_-{\pi} & (G/T)_T \ar[l]_-j \ar[d]^{\pi_{{}_T}} \\
(\pt)_G &  (\pt)_T  \ar[l]}
};
(2,-4)*{BG =};
(48,-4)*{= BT.};
\end{xy}

\bigskip
\noindent
By the push-pull formula
(\cite[Prop.~8.3]{borel--hirzebruch} or \cite[Lem.~1.5]{chern}), this diagram induces a 
commutative diagram in cohomology
\[
\xymatrix{
{H^*(BT)\ \ } \ar@{>->}[r]^-{j^*} \ar[d]_-{\pi_*} & H_T^*(G/T) \ar[d]^-{\pi_{{}_{T*}}} \\
{H^*(BG)\ } \ar@{^{(}->}[r] &H^*(BT),}
\]
where the horizontal maps are injections by the discussion of Section~\ref{s:relation}.
For $w \in W \subset G/T$, let $i_w\colon \{w\} \to G/T$ be
the inclusion map and $i_w^*\colon H_T^*(G/T) \to H_T^*\big(\{ w\}\big) = H^*(BT)$
the restriction map in equivariant cohomology.
For $x=b(u) \in H^*(BT)$, recall that $j^*b(u) = b(\ty)$ and $i_w^* b(\ty) = w\cdot b(u)$.
As in the proof of Theorem~\ref{t:pushforward}, 
by applying the equivariant localization theorem to the $T$-manifold $G/T$,
we obtain
\begin{align*}
\pi_*x &=\pi_*b(u) = \pi_{{}_{T*}}j^*b(u) = \pi_{{}_{T*}}b(\ty) 
=\sum_{w \in W} \frac{i_w^* b(\tilde{y})}{e^T(\nu_w)}\\
&= \sum_{w\in W}
w\cdot \left(
\dfrac{b(u)}{\prod_{\ga \in \tri^+} c_1(S_{\ga})} \right)
= \dfrac{\sum_{w\in W} (-1)^w w\cdot b(u)}
{\prod_{\ga \in \tri^+} c_1(S_{\ga})}
= \dfrac{\sum_{w\in W} (-1)^w w\cdot x}
{\prod_{\ga \in \tri^+} a_i}. \tag*{\qed}
\end{align*}
\epf

\subsection{The Associated Complete Flag Bundle} \label{ss:flag}

Suppose $V\rightarrow M$ is a $\cinf$ complex vector bundle of rank $n$.  
Let $f\colon  \Fl(V) \rightarrow M$ be the associated bundle of 
complete flags in the fibers of $V$.
It is a fiber bundle with fiber $G/T$, where $G$ is the unitary group $\Unit(n)$
and $T$ is the maximal torus
\[
T= \left\{
\left.
t= \begin{bmatrix}
t_1 & & \\
& \ddots & \\
& & t_n
\end{bmatrix}\ \right|
t_i \in \Unit(1) \right\}  = \Unit(1) \times \cdots \times \Unit(1) = \Unit(1)^n.
\]
The Weyl group of $T$ in $\Unit(n)$ is $S_n$, the symmetric group
on $n$ letters \cite[Th.\ IV.3.2, p.~170]{brocker--tomDieck}.

Consider the basis $\chi_1, \dots, \chi_n$ for the characters 
of $T$, where $\chi_i(t)= t_i$.  
A simple calculation of $tAt\inv$, where $t \in T$ and 
$A=[a_{ij}]$ is an $n \times n$ matrix, shows that the roots
of $\Unit(n)$ are $\chi_i \chi_j\inv$, $i \ne j$,
or in the additive notation of this paper, $\chi_i - \chi_j$.
(The root $\chi_i - \chi_j$ is the function$: T \to \Unit(1)$ given by
$t^{\chi_i - \chi_j} = \chi_i(t) \chi_j(t)\inv = t_i t_j\inv$.)
These are the \emph{global roots}, not the \emph{infinitesimal
roots}, of a Lie group \cite[Def.~V.1.3, p.~185]{brocker--tomDieck}.
A choice of positive roots for $\Unit(n)$ is
\[
\tri^+ = \{ \chi_i - \chi_j \mid 1 \le i < j \le n \}.
\]

Recall from \eqref{e:gtg} that $(G/T)_G = BT$.  By 
Lemma \ref{l:universal}, there are bundle maps $\bar{h}$ and $\uh$,
\smallskip
\[
\begin{xy}
(0,5)*{\xymatrix{
\Fl(V) \ar[d] \ar[r]^-{\bar{h}} & (G/T)_G \ar[d] \\
M \ar[r]_-{\uh} & BG,}
};
%(36,10)*{=BT};
(34,10)*{=BT};
\end{xy}
\] 
\vspace{.3in}

\noindent
and correspondingly, ring homomorphisms in cohomology
\vspace{.1in}
\[
\begin{xy}
(0,5)*{\xymatrix{
H^*\big(\Fl(V)\big) & H^*(BT) \ar[l]_(.4){\bar{h}^*} \\
H^*(M) \ar[u]^{f^*} & H^*(BG) \ar[l]^{\underline{h}^*} \ar@{^{(}->}+<0ex,+2.5ex>;[u]}
};
%(49.25,10)*{\simeq \Q[u_1, \ldots, u_n]};
%(51.5,-4.5)*{\simeq \Q[u_1, \ldots, u_n]^{S_n}.};
(45,10)*{\simeq \Q[u_1, \ldots, u_n]};
(47,-4.5)*{\simeq \Q[u_1, \ldots, u_n]^{S_n}.};
\end{xy}
\]
\vspace{.1in}

\noindent
By \eqref{e:cohombg}, the vertical map on the right is an injection.
By Theorem~\ref{t:formal}, the elements $a_i := \bar{h}^*(u_i)
\in H^2\big(\Fl(V)\big)$ generate $H^*\big(\Fl(V)\big)$ as an algebra over $H^*(M)$.

We will now deduce from Theorem~\ref{t:gysingt} a formula for the 
pushforward map $f_*$.

\prop 
For the associated complete flag bundle $f\colon \Fl(V) \to M$,
if $b(u) \in H^*(BT) = \Q[u_1, \ldots, u_{\ell}]$, then
\[
f^* f_* b(a) = \sum_{w\in S_n} w\cdot
\left(
\dfrac{b(a)}
{\prod_{i < j} (a_i - a_j)} \right),
\]
where $w\cdot b(a_1, \ldots, a_n) = b\big(a_{w(1)}, \ldots, a_{w(n)}\big)$.
\eprop

\pf 
Since $\bar{h}^*\colon H^*(BT) \to H^*\big(\Fl(V)\big)$ is a ring homomorphism,
for $b(u) \in H^*(BT) =\Q[u_1, \ldots, u_n]$,
\[
b(a) = b(a_1, \dots, a_n) = b(\oh^*u_1, \dots, \oh^*u_n) = 
\oh^*b(u_1, \dots, u_n) = \oh^*b(u) \in H^*\big(\Fl(V)\big).
\]
By Theorem~\ref{t:gysingt},
\begin{align*}
f^*f_* b(a) &= f^*f_*\bar{h}^* b(u)\\
&= \bar{h}^* \sum_{w\in S_n} w\cdot
\left(
\frac{b(u)}
{\prod_{\ga\in \tri^+ } c_1 (S_{\ga})}  \right).
\end{align*}
Since $\bar{h}^*$ commutes with $w$ (p.~\pageref{p:whomomorphism}),
\[
\bar{h}^* \big(w \cdot b(u)\big) = w\cdot \big(\bar{h}^*b(u)\big) = w\cdot b(a).
\]
If $\ga \in \tri^+$, then $\ga = \chi_i - \chi_j$
for some $1 \le i < j \le n$,
so that
\[
\bar{h}^*c_1(S_{\ga}) = \bar{h}^* c_1( S_{\chi_i - \chi_j}) = \bar{h}^*(u_i - u_j) = a_i - a_j.
\]
%and
%\[
%\bar{h}^* \left( w \cdot 
%\prod_{\ga\in \tri^+} c_1 (S_{\ga}) \right) 
%= w \cdot \prod_{i< j} (a_i - a_j).
%\]
Hence,
\[
f^*f_*b(a) = \sum_{w\in S_n} w\cdot \left(
\dfrac{b(a)}
{\prod_{i<j} (a_i - a_j)} \right),
\]
which agrees with \cite[Section 4.1, p.~41]{fulton--pragacz}.
\qed\epf

\subsection{The Associated Grassmann Bundle}

For a complex vector bundle $V\rightarrow M$ of rank $n$, 
the associated Grassmann bundle $f\colon  G(k,V) \rightarrow M$ of
$k$-planes in the fibers of $V$ is a fiber 
bundle with fiber the Grassmannian $G(k,n) = G/H$, where 
\[
G=\Unit(n) \quad\text{and}\quad H = \Unit(k) \times \Unit(n-k).
\]
A maximal torus contained in $H$ is $T = \Unit(1)^n$.
The Weyl groups of $T$ in $G$ and $H$ are 
\[
W_G = S_n \quad\text{and} \quad W_H = S_k \times S_{n-k}.
\]

If we let $\chi_i$ and $\tri^+$ be as in Subsection~\ref{ss:flag},
a choice of positive roots for $T$ in the subgroup $H$ is
\[
\tri^+(H) = \{ \chi_i - \chi_j \mid 1 \le i < j \le k \} \cup
\{ \chi_i - \chi_j \mid k+1 \le i < j \le n \}.
\]
By \eqref{e:cohombg},
\[
H^*(BH) = H^*(BT)^{W_H} = \Q [u_1, \dots, u_n]^{S_k \times 
S_{n-k}},
\]
where $u_i = c_1(S_{\chi_i})$.

Over $G(k,V)$ there are a tautological subbundle $S$ and a 
tautological quotient bundle $Q$, with total Chern classes
\begin{align*}
c(S) &= 1 + c_1(S) + \cdots + c_k(S) = \prod_{i=1}^k (1+a_i),\\
c(Q) &= 1+ c_1(Q) + \cdots + c_{n-k}(Q) = \prod_{i=k+1}^n (1+a_i).
\end{align*}
The $a_i$ for $1 \le i \le k$ are called the \term{Chern roots} of $S$,
and the $a_i$ for $k+1 \le i \le n$ the \term{Chern roots} of $Q$
(see \cite[\S 21, The Splitting Principle]{bott--tu}).
These $a_i$ are not cohomology classes on $G(k,V)$,
but are classes on the complete flag bundle $\Fl(V)$,
each of degree $2$.
The cohomology ring of $G(k,V)$ is
\begin{align*}
H^*\big(G(k,V)\big) &= \frac{H^*(M)\big[c_1(S), \ldots, c_k(S),
c_1(Q), \ldots, c_{n-k}(Q)\big]}{\big(c(s) c(Q) - f^*c(V)\big)}\\
&= \frac{H^*(M) \otimes_{\Q} \Q[ a_1, \dots, a_n]
^{S_k \times S_{n-k}}}
{\big(\prod (1+a_i) - (1+e_1+ \dots + e_n)\big)},
\end{align*}
where $e_i$ is the $i$-th Chern class $c_i(f^*V)$.

\prop [\cite{pragacz}, Lemma 2.5]
For the associated Grassmann bundle $f\colon  G(k,V) \rightarrow M$,
if $b(a)=b(a_1, \dots, a_n) \in \Q[a_1, \ldots, a_n]^{S_k \times S_{n-k}}$,
then its pushforward under $f_*$ is given by
\[
f^*f_*b(a) = \sum_{w \in S_n /(S_k \times S_{n-k})} w\cdot \left(
\dfrac{b(a)}
{\prod_{i=1}^k \prod_{j=k+1}^n (a_i - a_j)} \right).
\]
\eprop

\pf 
For $p \in M$, denote the fiber of the vector bundle $V$ over $p$ by $V_p$
and let $\Fl(V) \to G(k,V)$ be the natural map that sends
a complete flag $\Lambda_1 \subset \cdots \subset \Lambda_n = V_p$
where
$\dim_{\C} \Lambda_i = i$ to the partial flag $\Lambda_k \subset V_p$.
Let $G = \Unit(n)$, $H = \Unit(k) \times \Unit(n-k)$, and $T = \Unit(1)^n$.
This map $\Fl(V) \to G(k,V)$ is a fiber bundle with fiber $H/T$ and
group $G$.
If $P$ is the bundle of all unitary frames of the vector bundle $V$,
then $P$ is a principal $G$-bundle, and $\Fl(V) = P \times_G (G/T)$
and $G(k,V) = P \times_G (G/H)$ are the associated fiber bundles
with fiber $G/T$ and $G/H$ respectively.

Recall that $(G/H)_G = EG \times_G (G/H) \simeq BH$.  
As in Lemma~\ref{l:universal}, the classifying map
$\underline{h}\colon M \to BG$ of the principal bundle $P \to M$
induces a commutative diagram of bundle maps
\smallskip
\[
\begin{xy}
(0,5)*{\xymatrix{
\Fl(V) \ar[d] \ar[r]^{\bar{h}} & (G/T)_G \ar[d] \\
G(k,V) \ar[r]^-h \ar[d]_-f & (G/H)_G \ar[d]^-{\pi} \\
M \ar[r]_-{\underline{h}} & BG,}
};
%(38,10)*{\simeq BT};
%(38,-4)*{\simeq BH};
(36,10)*{\simeq BT};
(36,-4)*{\simeq BH};
\end{xy}
\] 
\vspace{.6in}

\noindent
and correspondingly, a diagram of ring homomorphisms in cohomology
\vspace{.1in}
\[
\begin{xy}
(0,5)*{\xymatrix{
H^*\big(\Fl(V)\big) & H^*(BT) \ar[l]_(.4){\bar{h}^*} \\
H^*\big(G(k,V)\big) \ar[u] & H^*(BH) \ar[l]_(.4){h^*}  \ar@{}[u]|{\bigcup}\\
H^*(M) \ar[u]^{f^*} & H^*(BG) \ar[l]^{\underline{h}^*} \ar@{}[u]|{\bigcup}}
};
%(50.5,10)*{\simeq \Q[u_1, \ldots, u_n]};
%(57.5,-5)*{\simeq \Q[u_1, \ldots, u_n]^{S_k \times S_{n-k}}};
%(53,-20)*{\simeq \Q[u_1, \ldots, u_n]^{S_n}.};
(46,10)*{\simeq \Q[u_1, \ldots, u_n]};
(51,-5)*{\simeq \Q[u_1, \ldots, u_n]^{S_k \times S_{n-k}}};
(48,-20)*{\simeq \Q[u_1, \ldots, u_n]^{S_n}.};
\end{xy}
\]
\vspace{.1in}

\noindent
As in Subsection~\ref{ss:flag}, the Chern roots $a_i$ are precisely $\bar{h}^*(u_i)$.

By Theorem~\ref{t:gysingh},
\begin{align*}
f^*f_* b(a) &= f^*f_*h^* b(u)\\
&= h^* \sum_{w\in W_G/W_H} w\cdot
\left(
\frac{b(u)}
{\prod_{\ga\in \tri^+ - \tri^+(H)} c_1 (S_{\ga})}  \right)\\
&= \sum_{w\in S_n/(S_k \times S_{n-k})} w\cdot \left(
\dfrac{b(a)}
{\prod_{i=1}^k \prod_{j=k+1}^n (a_i - a_j)} \right). \tag*{\qed}
\end{align*}
\epf

% \section{Transporting pushforward formulas to equivariant 
% cohomology} \label{s:transport}
% 
% \begin{center}
% To be written
% \end{center}

\section{Symmetrizing Operators}

Interpolation theory is concerned with questions such as how to find
a polynomial on $\R^n$ with given values at finitely many given points.
In interpolation theory there are symmetrizing operators that
take a polynomial with certain symmetries to another polynomial
with a larger set of symmetries.
For example, the Lagrange--Sylvester symmetrizer takes
a polynomial symmetric in two sets of variables $x_1, \ldots, x_k$
and $x_{k+1}, \ldots, x_n$ separately to a polynomial
symmetric in all the variables $x_1, \ldots, x_n$.
A curious byproduct of our Theorem~\ref{t:pushforward} 
is that it provides a geometric 
interpretation and consequently a generalization of some 
symmetrizing operators in interpolation theory \cite{lascoux}.

Let $X_n = (x_1, \dots, x_n)$ be a sequence of variables and 
$\Z[X_n] = \Z[x_1, \dots, x_n]$
the polynomial ring over $\Z$ generated by
$x_1, \ldots, x_n$.   
The \term{Lagrange--Sylvester 
symmetrizer} is the operator 
$\Delta\colon \Z[X_n]^{S_k \times S_{n-k}} \rightarrow 
\Z[X_n]^{S_n}$ taking $b(x) \in 
\Z[X_n]^{S_k \times S_{n-k}}$ to
\[
\Delta b(x) = \sum_{w \in S_n / (S_k \times S_{n-k})} w \left( 
\frac{b(x)}
{\prod_{i=1}^k \prod_{j=k+1}^n (x_j -x_i)}\right).
\]

The \term{Jacobi symmetrizer} is the 
operator $\partial\colon  \Z[X_n] \rightarrow \Z[X_n]^{S_n}$
taking $b(x) \in \Z[X_n]$ to
\[
\partial b(x) = \sum_{w \in S_n} w \left(
\frac{b(x)}
{\prod_{i<j} (x_j- x_i)} \right).
\]

Let $G$ be a compact Lie group of rank $n$ and $H$ a closed subgroup
containing a maximal torus $T$ of $G$.
Let $W_H$ and $W_G$ be the Weyl groups of $T$ in 
$H$ and in $G$ respectively.
Theorem \ref{t:pushforward} suggests that to every compact 
Lie group $G$ and closed subgroup $H$ of maximal rank, one can associate 
a symmetrizing operator on the polynomial ring 
$\Z[X_n]^{W_H}$ as follows.

The map $\pi\colon  G/H \rightarrow \pt$ induces a pushforward map in 
$G$-equivariant cohomology,
\[
\pi_*\colon  H_G^*(G/H) \rightarrow H_G^*(\pt).
\]  
Now the $G$-equivariant cohomology with integer coefficients of $G/H$ is
\[
H_G^*(G/H) = H^*(BH) = \Z [X_n]^{W_H}
\]
and the $G$-equivariant cohomology with integer coefficients of a point is
\[
H_G^*(\pt) = H^*(BG) = \Z [ X_n]^{W_G}.
\]
For the action of $T$ on $G/H$, the fixed point set is 
$W_G/W_H$.
Let $\tri^+(H)$ be a set of positive roots of $H$, and $\tri^+$ a 
set of positive roots of $G$ containing $\tri^+(H)$. 
As in Section~\ref{s:complete}, $c_{ET/T}$ is the 
characteristic map of $BT = ET / T$.   
The equivariant Euler 
class of the normal bundle at the identity element of $G/H$ 
is $\prod_{\ga \in \tri^+- \tri^+(H)} 
c_{ET/T}(\ga) \in H^*(BT) \simeq \Z [ X_n]$
(see \cite[Th.~19]{tu}).
Following our computation for the Gysin maps of a complete
flag bundle and a Grassmann bundle (but with integer instead
of rational coefficients),
we define the 
symmetrizing operator
\[
\Box \colon 
 \Z[X_n]^{W_H} \rightarrow \Z[X_n]^{W_G}
\]
 to be the operator taking $b(x) \in \Z[X_n]^{W_H}$ to
 \[
 \Box\, 
b(x) = \sum_{w \in W_G/W_H} w 
 \left(
 \frac{b(x)}
 {\prod_{\ga \in \tri^+ - \tri^+(H)} c_{ET/T}(\ga)}\right).
 \]
The Lagrange--Sylvester symmetrizer is the special case 
$G=\Unit(n)$, $H=\Unit(k) \times \Unit(n-k)$, and $T = \Unit(1)^n$,
and the Jacobi 
symmetrizer is the special case $G= \Unit(n)$ and $H= T = \Unit(1)^n$.
The meaning of the generalized symmetrizing operators awaits further
investigation.

% \section{Chow groups}
% 
% The preceding discussion is for the real equivariant 
% cohomology in the $\cinf$ category.  Edidin, Graham, and 
% Totaro have developed a parallel theory of equivariant Chow 
% groups  for algebraic varieties, including an analogue for 
% equivariant Chow groups of the localization formula.  Since 
% the formalisms in the algebraic and $\cinf$ categories are 
% the same, all of our results go through for the pushforward 
% in the Chow groups of a fiber bundle of varieties.  We refer 
% to \cite{graham} for details of this transition. 
% 

\end{document}